\documentstyle[11pt,leqno, pb-diagram]{article}
\input amssym.def
\catcode`\û=\active
\catcode`\…=\active
\catcode`\Ù=\active
\catcode`\–=\active
\catcode`\ð=\active
\catcode`\š=\active
\catcode`\ö=\active
\catcode`\û=\active
\letû=\ss
\def…{\"A}
\defÙ{\"a}
\def–{\"Oe}
\defð{\"o}
\defš{\"U}
\defö{\"u}
\title{Cartan Spinor Bundles on Manifolds.  \footnote {Supported by the SFB 288 of the DFG.}}
\author{Thomas Friedrich, Berlin}
\date{\today}

\begin{document}

\newcommand{\rpn}{{\Bbb R}{\Bbb P}^n}

\maketitle


\section{Introduction.}

Spinor fields and Dirac operators on Riemannian manifolds $M^n$ can be introduced by means of a reduction of the frame bundle from the structure group $O(n)$ to one of the groups $Spin(n)$, $Pin(n)$ or $Spin^{\Bbb C} (n)$. In any case the existence of a corresponding reduction imposes topological restrictions on the manifold $M^n$. The aim of the present paper is the construction of spinor bundles of Cartan type over certain non-orientable manifolds. The bundles under consideration do not split into subbundles invariant under the action of the Clifford algebra  and therefore they are not induced from $Pin^{\Bbb C}$ structure of the manifold. Moreover, we study the case of the real projective space $\rpn$ and its spinor bundle in more detail. These bundles over $\rpn$ admit a suitable metric connection $\nabla^S$ and  the corresponding Dirac operator. In particular it turns out  that there are Killing spinors in the constructed twisted spinor bundle over $\rpn$ for all dimensions.\\

Let us fix some notations. In case $n \equiv 0 \, mod \, 2$ we denote by $\gamma : \mbox{Cliff}^{\Bbb C} (n) \stackrel{\sim}{\to} \mbox{End} (\Delta_n)$ the isomorphism between the Clifford algebra and the algebra of endomorphisms of the space $\Delta_n$ of all Dirac spinor. In this way we obtain the so called {\it Dirac representation} $\gamma$ of the Clifford algebra. In case $n \equiv 1 \, mod \, 2$ we have the {\it Pauli representation} $\gamma : \mbox{Cliff}^{\Bbb C} (n) \to \mbox{End} (\Delta_n)$. Denote by $\alpha : \mbox{Cliff}^{\Bbb C} (n) \to \mbox{Cliff}^{\Bbb C} (n)$ the canonical involution. Then

\[ \gamma \oplus (\gamma \circ \alpha) : \mbox{Cliff}^{\Bbb C} (n) \stackrel{\sim}{\to} \mbox{End} (\Delta_n) \oplus \mbox{End} (\Delta_n) \]

is an isomorphism and the corresponding representation of the Clifford algebra in $\Delta_n \oplus \Delta_n$ is the so called {\it Cartan representation} (see [Fr2], [Tr]).\\

\section{The construction of spinor bundles of Cartan type.}

We consider {\it odd-dimensional} Riemannian manifolds $M^{2k+1}$ only. A spinor bundle of Cartan type is a $2^{k+1}$-dimensional complex Clifford bundle $S$ realizing at any point the Cartan representation. $\gamma \oplus (\gamma \circ \alpha): \mbox{Cliff}^{\Bbb C} ({\Bbb R}^{2k+1}) \stackrel{\sim}{\to} \mbox{End} (\Delta_{2k+1})  \oplus \mbox{End} (\Delta_{2k+1})$ of the Clifford algebra (see [Tr]). Any spinor bundle of Cartan type over an orientable manifold splits into two Clifford subbundles. Indeed, the volume form $dM^{2k+1} = e_1 \cdot \ldots \cdot e_{2k+1} \in \mbox{Cliff} (T(M^{2k+1}))$ commutes with the Clifford multiplication. Since

\[ \alpha(e_1 \cdot \ldots \cdot e_{2k+1})= -e_1 \cdot \ldots \cdot e_{2k+1} \]

$dM^{2k+1}$ acts on $S$ with two different eigenvalues. Therefore $dM^{2k+1}$ defines a splitting of the bundle $S$ into two Clifford subbundles.\\

Consider a simply-connected, odd-dimensional Riemannian spin manifold $M^{2k+1}$ and denote by $S$ its spinor bundle. Any isometry $\gamma : M^{2k+1} \to M^{2k+1}$ admits two lifts $\hat{\gamma}_{\pm} :S \to S$ into the spinor bundle $S$. We once again describe the construction of $\hat{\gamma}_{\pm}$. The differential $d \gamma $ (resp. $(-d \gamma) - M$ is odd-dimensional !) maps the frame bundle into itself if $\gamma$ preserves the orientation  (resp. does not preserve the orientation). Moreover, $d \gamma $ (or $(-d \gamma)$) lifts into the spin structure of $M^{2k+1}$ and defines two lifts $\hat{\gamma}_{\pm}$. The lifts $\hat{\gamma}_{\pm}$ commute (anticommute) with the Clifford multiplication, i.e.

\[ \hat{\gamma}_{\pm} (X \cdot \varphi) = \deg (\gamma ) d \gamma (X) \cdot \hat{\gamma}_{\pm} (\varphi)  \]

(see [Fr1]).\\

Let $\Gamma$ be a discrete subgroup of isometries and suppose that $\Gamma$ acts freely on $M^{2k+1}$. We consider the spinor bundle of Cartan type 

\[ S^* = S \oplus S \]

over $M^{2k+1}$. The Clifford multiplication  in $S^*$ is given by the formula

\[ X \cdot \left( \begin{array}{c} \varphi\\ \psi \end{array} \right) = \left( \begin{array}{c} X \cdot \varphi\\ -X \cdot \psi \end{array} \right) \hspace{2cm} X \in T(M^{2k+1}) \]

where $X \cdot \varphi$ denotes the Clifford multiplication in the spinor bundle $S$. Any isometry $\gamma \in \Gamma$ admits the following 4 lifts into the bundle $S^*$:

\[ \left( \begin{array}{cc} \hat{\gamma}_{\pm} & 0\\ 0 & \hat{\gamma}_{\pm} \end{array} \right) \hspace{2cm} \mbox{in case $\gamma$ preserves the orientation.} \]

\[ \left( \begin{array}{cc} 0 & \hat{\gamma}_{\pm}\\  \hat{\gamma}_{\pm} & 0\end{array} \right) \hspace{2cm} \mbox{in case $\gamma$ reverses the orientation.} \]

Any of these lifts commutes with the Clifford multiplication in the bundle $S^*$. \\

We introduce now a family of automorphisms $\Pi (\Theta_1, \Theta_2, \varepsilon_1, \varepsilon_2, \gamma)$ of the bundle $S^*$ depending on four parameters $0 \le \Theta_1 \le 2 \pi$, $0 \le \Theta_2 \le 2 \pi$ and $\varepsilon_1, \varepsilon_2 \in \{ \pm1 \}$. If $\gamma$ preserves the orientation of the manifold we define

\[ \Pi (\Theta_1, \Theta_2, \varepsilon_1, \varepsilon_2, \gamma)= \left( \begin{array}{cc} e^{i \Theta_1} \hat{\gamma}_{\varepsilon_1} & 0 \\ \mbox{}\\
0 & e^{i \Theta_2} \hat{\gamma}_{\varepsilon_2} \end{array} \right) , \]

otherwise let us introduce the automorphism

\[ \Pi (\Theta_1, \Theta_2, \varepsilon_1, \varepsilon_2, \gamma)= \left( \begin{array}{cc} 0 & e^{i \Theta_2} \hat{\gamma}_{\varepsilon_2}  \\ \mbox{}\\
 e^{i \Theta_1} \hat{\gamma}_{\varepsilon_1} & 0  \end{array} \right) . \]

Suppose that $\varepsilon : \Gamma \to \mbox{Aut} (S^*)$ is a homomorphism such that for any $\gamma \in \Gamma$ the automorphism $\varepsilon (\gamma)$ coincides with one of the lifts $\Pi (\Theta_1, \Theta_2, \varepsilon_1, \varepsilon_2, \gamma)$ of $\gamma$ into $S^*$. Then

\[\bar{S} =S^*/ \varepsilon (\Gamma) \]

is a Clifford bundle of Cartan type over the Riemannian manifold $\bar{M}^{2k+1} = M^{2k+1} / \Gamma$. The bundle $\bar{S}$ splits into two parts invariant under the action of the Clifford bundle $\mbox{Cliff} (T(M^{2k+1}))$ if and only if $\bar{M}^{2k+1}$ is orientable. For a non-orientable manifold $\bar{M}^{2k+1}$ the non-splitting bundle $\bar{S}$ is not induced from a $Pin^{\Bbb C} (2k+1)$ structure of $\bar{M}^{2k+1}$.\\

{\bf Example:} $\Gamma = {\Bbb Z}_p$.\\

Suppose $\Gamma$ is generated by one involution $\gamma : M^{2k+1} \to M^{2k+1}$ with the lifts $\hat{\gamma}_{\pm}$. If $\gamma$ preserves the orientation and $\hat{\gamma}^2_{\pm} = \mbox{Id}$ the manifold $\bar{M}^{2k+1} = M^{2k+1} / {\Bbb Z}_2$ admits a spin structure. In case $\hat{\gamma}^2_{\pm} = \mbox{-Id}$ \, \, $\bar{M}^{2k+1}$ is orientable, but does not admit a spin structure. Nevertheless we can construct a spinor bundle using the homomorphism

\[ \varepsilon (\gamma) = \left( \begin{array}{cc} i \hat{\gamma}_+ & 0 \\ 0 & i \hat{\gamma}_+ \end{array} \right) . \]

If $\gamma$ does not preserve the orientation we define the homomorphism $\varepsilon : {\Bbb Z}_2 \to \mbox{Aut} (S^*)$ by the formula

\[ \varepsilon (\gamma) = \left( \begin{array}{cc} 0 & \hat{\gamma}_+ \\ \hat{\gamma}_+ &0 \end{array} \right) \hspace{2cm} \mbox{if $\hat{\gamma}^2_+ = \mbox{Id}$} \]

\[ \varepsilon (\gamma) = \left( \begin{array}{cc} 0 & i \hat{\gamma}_+ \\ i  \hat{\gamma}_+ &0 \end{array} \right) \hspace{2cm} \mbox{if $\hat{\gamma}^2_+ = \mbox{-Id}$} .  \]

A similar construction for an arbitrary cyclic group ${\Bbb Z}_p$ yields the result: \\

{\it Let $\bar{M}^{2k+1}$ be a Riemannian manifold with fundamental group $\pi_1 (\bar{M}^{2k+1}) = {\Bbb Z}_p$ and suppose that the universal covering $M^{2k+1}$ admits a spin structure. Then over $\bar{M}^{2k+1}$ there exists a spinor bundle $\bar{S}$ of Cartan type. The bundle $\bar{S}$ splits into two subbundles invariant under the action of the Clifford bundle} $\mbox{Cliff} (T(\bar{M}^{2k+1}))$ {\it if and only if $\bar{M}^{2k+1}$ is  orientable.}\\

\newfont{\graf}{eufm10}
\newcommand{\D}{\mbox{\graf D}}

Let $D: \Gamma (S) \to \Gamma (S)$ be the Dirac operator acting on sections of the spinor bundle $S$ over $M^{2k+1}$. The Dirac operator $\D$ of the spinor bundle $S^*$ of Cartan type is given by

\[ \D = \left( \begin{array}{cc}D &0\\0 & - D \end{array} \right) . \]

Therefore the eigenspaces $E_{\lambda} (\D)$ consists of pairs $\left( \begin{array}{c} \varphi\\ \psi \end{array} \right)$ of eigenspinors of $D$, i.e.

\[ E_{\lambda} (\D) =E_{\lambda} (D) \oplus E_{- \lambda} (D) . \]

The group $\Gamma $ acts on $E_{\lambda} (\D)$ by the matrices $\varepsilon (\gamma)$ and the eigenspace of the Dirac operator $\bar{\D}$ on $\bar{M}^{2k+1}$ coincides with the subspace of $\varepsilon (\Gamma)$-invariant pairs in $E_{\lambda} (\D)$. In particular, in case $\Gamma = {\Bbb Z}_2$ and $\bar{M}^{2k+1}$ is non-orientable we have only one transformation

\[€\varepsilon (\gamma) = \left( \begin{array}{cc} 0 & \hat{\gamma}_+ \\ \hat{\gamma}_+ & 0 \end{array} \right) \quad \mbox{or} \quad €\varepsilon (\gamma) = \left( \begin{array}{cc} 0 & i \hat{\gamma}_+ \\ i \hat{\gamma}_+ & 0 \end{array} \right) . \]

The map $E_{\lambda} (D) \ni \varphi \longmapsto \left( \begin{array}{c} \varphi\\ \hat{\gamma}_+ (\varphi) \end{array} \right) \in E_{\lambda} (\bar{\D})$ (resp. $\varphi \longmapsto \left( \begin{array}{c} \varphi\\ i \hat{\gamma}_+ (\varphi) \end{array} \right)$) defines an isomorphism between the corresponding eigenspaces. Therefore the spectrum of the operator $(\bar{\D}, \bar{S}, \bar{M}^{2k+1})$ coincides with the spectrum of the operator $(D, S; M^{2k+1})$.\\

We thus constructed non-splitting spinor bundles over non-orientable manifolds $\bar{M}^{2k+1}$ via coverings. We mention that other situations give similarly bundles with the described properties. For example, any immersion  of $\bar{M}^{2k+1}$ into an oriented spin manifold $N^{2k+2}$ also defines a spinor bundle of Cartan type over $\bar{M}^{2k+1}$. Indeed, let $f : \bar{M}^{2k+1} \to N^{2k+2}$ be the immersion and consider the Dirac spinor bundle $S$ of $N^{2k+2}$. Then the induced bundle $\bar{S} = f^* (S)$ inherits the structure of a spinor bundle of Cartan type over $\bar{M}^{2k+1}$.\\

\section{The Spinor Bundle over $\rpn$.}

We start with the $n$-dimensional sphere

\[ S^n = \{ x \in {\Bbb R}^{n+1} : |x| =1 \} \]

and we denote by $S^*$ the $2^{\left[ \frac{n+1}{2} \right]}$-dimensional bundle

\[ S^* = S^n \times \Delta_{n+1} . \]

The tangent bundle of the sphere

\[ T(S^n)= \{ (x,t) \in S^n \times {\Bbb R}^{n+1} : \langle x,t \rangle =0 \} \]

acts by Clifford multiplication

\[ \mu  ((x,t) \otimes (x, \varphi)):=(x,t \cdot \varphi) \]

on $S^*$. Therefore $S^*$ is a Clifford bundle over $S^n$. Let $g: S^n \to S^n, g(x)=-x$ be the antipodal map. Its differential $dg: T(S^n) \to T(S^n)$ acts on a tangent vector $(x,t) \in T(S^n)$ by 

\[ dg(x,t) =(-x,-t) . \]

Let us define a lift $\hat{g}$ of the antipodal map into the bundle $S^*$ by the formula

\[ \hat{g} (x, \varphi)=(-x, x \cdot \varphi) . \]

Then $\hat{g}$ is again an involution, $(\hat{g})^2 =\mbox{Id}$. Moreover, $\hat{g}$ respects the Clifford multiplication, i.e. the diagram

\[ 
\begin{diagram}
\node{T(S^n) \otimes S^*} \arrow{s,b}{(dg) \otimes \hat{g}} \arrow{e,t}{\mu} \node{S^*} \arrow{s,b}{\hat{g}}\\
\node{T(S^n) \otimes S^*} \arrow{e,t}{\mu} \node{S^*}
\end{diagram}
\]

commutes. Indeed, given $(x,t) \otimes (x, \varphi) \in T(S^n) \otimes S^*$ we have

\[ \mu \circ (dg \otimes \hat{g} )((x,t) \otimes (x, \varphi))= \mu ((-x,-t) \otimes (-x, x \cdot \varphi)) = (-x,-t \cdot x \cdot \varphi ) . \]

An the other hand,

\[ \hat{g} \mu ((x,t) \otimes (x, \varphi))= \hat{g} ((x, t \cdot \varphi))=(-x, x \cdot t \cdot \varphi) . \]

Since $x$ and $t$ are orthogonal vectors in ${\Bbb R}^{n+1}$ we have in the Clifford algebra the relation

\[ x \cdot t + t \cdot x =0 , \]

i.e. the mentioned diagram commutes.\\

Over the real projective space $\rpn = S^n/g$ we define the bundle

\[ S=S^*/ \hat{g} . \]

The tangent bundle $T(\rpn)$ can be identified with $T(\rpn)=T(S^n)/dg$ and the above discussed property of the lift $\hat{g} : S^* \to S^*$ implies that we obtain a well-defined Clifford multiplication

\[ \mu : T(\rpn) \otimes S \to S . \]

At any point of $\rpn$ the bundle $S$ realizes the Dirac (in case $n \equiv 0 \, mod \, 2$) or the Cartan (in case $n \equiv 1 \, mod \, 2$) representation of the $n$-dimensional Clifford algebra.\\

Let us consider the even-dimensional case, $n =2k$. The Dirac  representation of the Clifford algebra decomposes into two irreducible representations with respect to the action of the even part $\mbox{Cliff}^{\Bbb C}_o (n)$ of the Clifford algebra. Therefore the bundle $S$ decomposes at any fixed point as a Clifford module too. However, globally the bundle $S$ over $\rpn$ does not split as a $\mbox{Cliff}_o (T(\rpn))$-bundle. Indeed suppose that there exists a decomposition $S=S_1 \oplus S_2$ invariant under the action of the Clifford bundle. Then we obtain a corresponding decomposition $S^* =S^*_1 \oplus S^*_2$ of the bundle $S^*$ over $S^n$. The spin module $\Delta_{n+1}$ decomposes in a unique way as a $\mbox{Cliff}_o (T_x S^n)$-module

\[ \Delta_{n+1} = \Delta^+_{n+1} (x) \oplus \Delta^-_{n+1} (x) . \]

We describe this decomposition in an explicit way. Let $e_1, \ldots, e_n$ an orthonormal basis of the tangent space $T_x (S^n)$ and consider the element

\[ f = i^{n/2} \, \, e_1 \cdot \ldots \cdot e_n : \Delta_{n+1} \to \Delta_{n+1} . \]

Then $f^2 = \mbox{Id}$ and $\Delta_{n+1}$ decomposes into the eigenspaces of $f, \Delta_{n+1} = \Delta^+_{n+1} (x) \oplus \Delta^-_{n+1}$ with 

\[ \Delta^{\pm}_{n+1} (x)= \{ \varphi \in \Delta_{n+1} : \quad f (\varphi) = \pm \varphi \} . \]

Now the volume form $e_1 \cdot \ldots \cdot e_n \cdot x$ of ${\Bbb R}^{n+1}$ acts on $\Delta_{n+1}$ by multiplication

\[ e_1 \cdot \ldots \cdot e_n \cdot x_{|\Delta_{n+1}} = \alpha_{n+1} . \]

Then $f =i^{n/2} \, \, e_1 \cdot \ldots \cdot e_n = - i^{n/2} \, \, \alpha_{n+1} x$ and the subspaces $\Delta^{\pm}_{n+1} (x)$ are given by

\[ \Delta^{\pm}_{n+1} (x) = \{ \varphi \in \Delta_{n+1} : \, \, \pm x \varphi = i^{n/2} \, \alpha_{n+1} \varphi \} . \]

The bundles $S^*_1$ and $S^*_2$ coincide therefore with the subbundles $\Delta^{\pm}_{n+1} (x)$. Since $S^*_i (i=1,2)$ are induced bundles by the covering $S^n \to \rpn =S^n/g$ they are invariant under the lift $\hat{g} :S^* \to S^*$. However, the lift $\hat{g}$ defined by the Clifford multiplication does not preserve the decomposition $\Delta_{n+1} = \Delta^+_{n+1} (x) \oplus \Delta^-_{n+1} (x)$, a contradiction.\\

We discuss now the odd-dimensional case, $n=2k+1$. Again we start with an orthonormal basis $e_1 , \ldots, e_n \in T_x(S^n)$ and we introduce the automorphism

\[ f= e_1 \cdot \ldots \cdot e_n : \Delta_{n+1} \to \Delta_{n+1} . \]

Then we have

\[ f^2 = (-1)^{k+1} \quad , \quad f \cdot x = - x \cdot f \]

and 

\[ f \cdot t = t \cdot f \quad \mbox{for all tangent vectors $t \in T_x (S^n)$} . \]

We decompose $\Delta_{n+1}$ into the eigenspaces of $f$

\[ \Delta_{n+1} = \Delta^+_{n+1} (x) \oplus \Delta^-_{n+1} (x) \quad , \quad \Delta^{\pm}_{n+1} (x) = \{ \varphi \in \Delta_{n+1} : e_1 \cdot \ldots \cdot e_n \varphi = \pm i^{k+1} \varphi \} . \]

Since $g$ commutes with all tangent vectors $t \in T_x (S^n)$ the subspaces $\Delta^{\pm}_{n+1} (x)$ are invariant under the action of the Clifford algebra $\mbox{Cliff} (T_x(S^n))$. The lift $\hat{g}$ preserves now $(n=2k+1)$ the decomposition $\Delta_{n+1}= \Delta^+_{n+1} (x) \oplus \Delta^-_{n+1} (x)$, i.e. the Clifford bundle $S$ over $\rpn$ decomposes into two Clifford subbundles.\\


Vector fields $V$ on the projective space $\rpn$ can be identified with maps $V:S \to {\Bbb R}^{n+1}$ such that

\[ V(-x)=-V(x) \quad , \quad \langle V(x), x \rangle =0 \quad x \in S^n . \]

The bundle $S$ over $\rpn$ arises from $S^* =S^n \times \Delta_{n+1}$ and the identification $\hat{g} (x, \varphi)=$ $=(-x, x \cdot \varphi)$. Therefore a section $\Phi \in \Gamma (S; \rpn)$ is a map $\Phi : S^n \to \Delta_{n+1}$ with the property\\

\mbox{} \hfill $ \Phi (-x)= x \cdot \Phi (x)$ .  \hfill $(*)$  \\

The formula

\[ (\nabla_V^S \Phi)(x)= d \Phi (V)(x) + \frac{1}{2} V(x) \cdot x \cdot \Phi (x) \]

defines a covariant derivative

\[ \nabla^S : \Gamma (T(\rpn)) \times \Gamma (S) \to \Gamma (S) . \]

Indeed, a simple calculation shows that $\nabla^S_V \Phi$ has the invariance property $(*)$ in case $\Phi$ and $V$ satisfy the corresponding transformation rules. Moreover, the covariant derivative $\nabla^S$ is compatible with the Clifford multiplication, i.e.

\[ \nabla^S_V (W \cdot \Phi)=(\nabla_VW) \cdot \Phi + W \cdot (\nabla^S_V \Phi) , \]

and $\nabla^S$ preserves the hermitian metric of the bundle $S$. To sum up, we see that the triple $(S, \nabla^S, \langle , \rangle)$ is a Dirac bundle over $\rpn$ (see [LM] for the general definition of a Dirac bundle). A direct calculation yields the formula

\[ R^S (V,W) \Phi = \frac{1}{2} \left(W \cdot V + \langle V,W \rangle \right) \cdot \Phi \]

for the curvature tensor $R^S$ of the connection $\nabla^S$. \\

{\bf Remark:} Let us compare the above defined covariant derivative $\nabla^S$ with the usual covariant derivative $\nabla$ in the usual spin bundle of a hypersurface in ${\Bbb R}^{n+1}$. The formula for $\nabla$ is

\[ \nabla_V \Phi = d \Phi (V) + \frac{1}{2} \mbox{II} \, (V) \cdot \vec{N} \cdot \Phi \]

i.e. over the sphere $S^n$ the introduced Dirac bundle coincides with the usual spinor bundle of the sphere. \\

The Dirac operator $D= \sum\limits^n_{i=1} e_i \nabla^S_{e_i}$ acting on sections $\Gamma (S; \rpn)$ is related to the Laplace operator $\Delta= - \sum\limits^n_{i=1} \nabla^S_{e_i} \nabla^S_{e_i} - \sum\limits^n_{i=1} \mbox{div} \, (e_i) \nabla^S_{e_i}$ by the well-known formula

\[ D^2 = \Delta + R \]

where the endomorphism $R$ (in case of a general Dirac bundle) is given by

\[ R= \frac{1}{2} \sum\limits_{j,k} e_j e_k R^S (e_j, e_k) . \]

In our situation we simply obtain

\[ D^2 = \Delta + \frac{1}{4} n(n-1) = \Delta + \frac{\tau}{4} \]

where $\tau$ is the scalar curvature of $\rpn$. The bundle $S$ over $\rpn$  admits Killing spinors, i.e. spinor fields $\Phi$ satisfying the differential equation

\[ \nabla^S_V \Phi = \lambda V \cdot \Phi \]

for some $\lambda \in {\Bbb R}^1$. Indeed, suppose that $\Phi$ is a Killing spinor in $\Gamma (S; \rpn)$. Then we obtain

    \[ R^S (V,W) \Phi = 2 \lambda^2 (WV + \langle V,W \rangle ) \cdot \Phi \]

using the differential equation. On the other hand, the formula for the curvature tensor implies

\[ R^S (V,W) \Phi = \frac{1}{2} (WV + \langle V,W \rangle ) \cdot \Phi , \]

i.e. $\lambda = \pm \frac{1}{2}$. Let us consider the spinor field

\[ \Phi (x)=(1-x) \Phi_0 \]

where $\Phi_0 \in \Delta_{n+1}$ is constant. Then

\[ \Phi (-x) =(1+x) \Phi_0 =x (1-x) \Phi_0 =x \Phi (x) \]

and therefore $\Phi \in \Gamma (S; \rpn)$ is a \,  section in the spinor bundle over $\rpn$. We calculate the covaraint derivative $\nabla^S_V \Phi$:

\begin{eqnarray*}
(\nabla^S_V \Phi)(x)&=& -V (x) \cdot \Phi_0 + \frac{1}{2} V(x) \cdot x \cdot (1-x) \Phi_0 =\\
&=& -V (x) \cdot \Phi_0 + \frac{1}{2} V(x) \cdot (1+x) \Phi_0 =\\
&=& \frac{1}{2} V(x) \cdot (-1+x) \cdot \Phi_0 = - \frac{1}{2} V(x) \cdot \Phi(x) , 
\end{eqnarray*}

i.e. $\Phi (x)$ is a Killing spinor with Killing number $\lambda= - \frac{1}{2}.$\\

{\bf Remark:} Killing spinors with Killing number $\lambda = + \frac{1}{2}$ on $\rpn$ appear if we use the identification $\hat{g}_- :S^* \to S^*$ defined by the formula

\[ \hat{g} (x, \varphi)=( -x, -x \cdot \varphi) . \]

$\hat{g}_-$ defines in a similar way a Clifford bundle $S_- =S^*/ \hat{g}_-$ whose sections are maps $\Phi :S^n \to \Delta_{n+1}$ with the property

\[ \Phi (-x) = - x \Phi (x) . \]

The spinor fields

\[ \Phi (x) =(1+x) \Phi_0 \]

are Killing spinors in the bundle $S_-$ over $\rpn$ with Killing number $\lambda = + \frac{1}{2}$. \\

Killing spinors on compact Riemannian spin manifolds with positive scalar curvature correspond to eigenspinors of the Dirac operator related to the smallest eigenvalue. In case ${\Bbb R}{\Bbb P}^{4k+3}$ the first eigenvalue of the Dirac operator on usual spinors has been calculated $\left( \lambda_1 = \pm \frac{1}{2} \sqrt{\frac{n}{n-1} \tau} \right)$ and there are Killing spinors (see [F1]). More general, we can calculate the spectrum of the Dirac operator acting on sections $\Phi \in \Gamma (S; \rpn)$ using the computation of the eigenvalues of $D$ on the sphere.  Let $\Phi : S^n \to \Delta_{n+1}$ be a section in the bundle $S^*$. We decompose $\Phi$ into $\Phi = \Phi_+ + \Phi_-$ with

\[ \Phi_+ (x) = \frac{\Phi(x) - x \cdot \Phi(-x)}{2} \quad , \quad \Phi_- (x) = \frac{\Phi (x) + x \cdot \Phi (-x)}{2} . \]

The sections $\Gamma (S; \rpn)$ are described by the condition $\Phi_- \equiv 0$. Now we apply the calculation of the spectrum of the Dirac operator on spheres as well as the realization of eigenspinors by polynomials (see [Su], [Tr]). Imposing the additional condition $\Phi_- \equiv 0$ we obtain the spectrum of the Dirac operator acting on the sections $\Gamma (S; \rpn)$ for arbitrary $n$.

\vspace{2cm}
Thomas Friedrich\\
Humboldt-UniversitÙt zu Berlin\\
Institut för Reine Mathematik\\
Sitz: Ziegelstraûe 13a\\
D-10099 Berlin\\
e-mail: friedric@mathematik.hu-berlin.de

\end{document}